\documentclass[a4paper,12pt]{article}
\usepackage[nottoc]{tocbibind}
\usepackage[cp1251]{inputenc}
\usepackage[english]{babel}
\usepackage{amsfonts}
\usepackage{amsfonts,amssymb}
\usepackage{mathenv}
\usepackage{amsmath}
\usepackage{amsthm}
\usepackage{wrapfig}
\usepackage{multirow}
\usepackage{multicol,multirow}
\usepackage{xcolor}
\usepackage[dvips]{graphicx}
\usepackage[dvips]{epsfig}
\usepackage{todonotes}
\usepackage{upgreek}
\usepackage{cite}
\usepackage{lineno,multirow,rotating,amsmath,setspace}
\usepackage{color}
\usepackage{amssymb, bm}
\usepackage{caption}
\usepackage{multicol}
\usepackage{geometry}
\usepackage{todonotes}
\usepackage{upgreek}
\modulolinenumbers[5]
\usepackage{nomencl}

\usepackage{graphicx}
\usepackage{hyperref}
\usepackage{lmodern}
\usepackage{url}
\usepackage{amssymb,amsthm,amsfonts,amstext}
\usepackage{amsmath}
\usepackage{graphicx}
\usepackage{enumerate}
\usepackage{a4wide}
\usepackage{todonotes}
\usepackage[mathscr]{eucal}
\newtheorem{remark}{Remark}[section]
\newtheorem{theorem}{Theorem}[section]

\newtheorem{corollary}{Corollary}[section]

\newtheorem{example}{Example}[section]
\renewcommand{\epsilon}{\varepsilon}

\numberwithin{equation}{section}

\renewcommand{\leq}{\leqslant}
\renewcommand{\geq}{\geqslant}

\begin{document}

\title{A new generalized integral transform and applications}

\author{Mohamed Akel\\Math. Department, Faculty of Science,\\ South Valley University, Qena 83523, Egypt\\makel@sci.svu.edu.eg}

\date{}
\maketitle
%
%

%

\begin{abstract}
	In this work, we introduce a new generalized integral transform involving many potentially known or new transforms as special cases. Basic properties of the new integral transform, that investigated in this work, include the existence theorem, the scaling property, elimination property a Parseval-type identity, and inversion formula. The relationships of the new transform with well-known transforms are characterized by integral identities. The new transform is applied to solve certain  initial boundary value problems. Some illustrative examples are given. The results established in this work extend and generalize recently published results.
\end{abstract}





\section{Introduction}\label{sec1}

For many decades, the integral transforms play a precious role in solving many  differential and integral equations. Using an appropriate integral transform helps to reduce differential and integral operators, from a considered domain into multiplication operators in another domain. Solving the deduced problem in the new domain, and then applying the inverse transform serve to invert the manipulated solution back to the required solution of the problem in its original domain (see, \cite{IT1,IT2,IT3,IT4,IT5,IT6,IT7,IT8,IT9,IT10,IT11,IT12,IT13}). 

The classical integral transforms used in solving differential equations, integral equations, and in analysis and the theory of functions are the Laplace transform, the Fourier integral transform, the Mellin transform \cite{IT8,IT9}. Besides, in the mathematical literature, there are many Laplace-type integral transforms such as the Laplace-Carson transform which is used in the railway engineering \cite{IT14}, the z-transform can be applied in signal processing \cite{IT15}, the Sumudu transform is used in engineering and many real-life problems \cite{IT10,IT16,0IT16}, the Hankel's and Weierstrass transform has been applied in heat and diffusion equations \cite{IT17,IT18}. In addition, we have the natural transform \cite{IT19,IT20,0IT20} and Yang transform \cite{IT13,IT21} used in many fields of physical science and engineering.

In this work, we study the following integral transform
\begin{align}\label{newIT}
\mathcal{M}_{\rho,m}[f(x)](u,v,\omega)=\int_0^\infty \frac{e^{-u x-v/x}}{\left(x^m +\omega^m\right)^\rho} f(\omega x) dx,
\end{align}
with $\rho\in\mathbb C,\,\operatorname{Re}\rho>0$, $m\in\mathbb N$ and $u,v\in\mathbb C$, and $\omega\in\mathbb R^+$ are the trasform variables.

Basic properties and applications of the new integral transform \eqref{newIT} are given in this study. Indeed, dualities of the new integral transform \eqref{newIT} with well-known integral transforms can be indicated as the following.
\begin{align}\label{newIT_SriIT}
\mathcal{M}_{\rho,m}[f(x)](u,0,\omega)=\mathbb{M}_{\rho,m}[f(x)](u,\omega),
\end{align}
where $\mathbb{M}_{\rho,m}[f(x)](u,\omega)$ is the Srivastava-Luo-Raina transform, defined in \cite{IT26} as
\begin{align}\label{SriIT}
\mathbb{M}_{\rho,m}[f(x)](u,\omega)=\int_0^\infty \frac{e^{-u x}}{\left(x^m +\omega^m\right)^\rho} f(\omega x) dx,
\end{align}
with $\rho\in\mathbb C,\,\operatorname{Re}\rho>0$, $m\in\mathbb N$, $u\in\mathbb C$, and $\omega\in\mathbb R^+$.

Setting $\rho=v=0$ in \eqref{newIT}, one gets the natural transform;
\begin{align}\label{naturalIT}
\mathcal{N}[f(x)](u,\omega)=\int_0^\infty e^{-u x}f(\omega x)dx,
\end{align}
whose properties were discussed in \cite{IT20} and its applications were given in \cite{IT22}. 
Thus, we have the following 
duality relations
$$\mathcal{N}[f(x)](u,\omega)=\mathcal{M}_{0,m}[f(x)](u,0,\omega),\; u,\omega>0,$$
\[
\mathcal{M}_{\rho,m}[f(x)](u,v,\omega)= \mathcal{N}\left[\frac{e^{-v\omega/x}f(x)}{\left(\frac{x^{m}}{\omega^{m}}+\omega^{m}\right)^{\rho}}\right](u,\omega),\; u,v, \omega>0,
\]
and,
\begin{align}\label{elimintateToNaturalIT}
\mathcal{M}_{\rho,m}\left[\left(\frac{x^{m}}{\omega^{m}}+\omega^{m}\right)^{\rho} f(x)\right](u,v,\omega)=
\mathcal{N}\left[e^{-v\omega/x}f(x)\right](u,\omega),\; u,v, \omega>0.
\end{align}
When $\rho=v=0$ and $\omega=1$, \eqref{newIT} reduces to the Laplace transform
\begin{align}\label{Laplace_T}
\mathcal{L}[f(x)](u)=\int_0^\infty e^{-ux}f(x)dx,\, \operatorname{Re}u>0.
\end{align}
So, from (\ref{newIT}) and (\ref{Laplace_T}) we have the following duality relations
\begin{align}\label{M_L duality}
\mathcal{L}[f(x)](u)&=\mathcal{M}_{0,m}[f(x)](u,0,1),\; \operatorname{Re}u>0,\nonumber\\
\mathcal{M}_{\rho,m}[f(x)](u,v,\omega)&=\mathcal{L}\left[\frac{e^{-v/x}f(\omega x)} {(x^{m}+\omega^{m})^{\rho}}\right](u),\; u,v,\omega>0,\\
\mathcal{M}_{\rho,m}[f(x)](u,v,\omega)&=\frac{1}{\omega}\mathcal{L}\left[\frac{e^{-v\omega/x}f(x)}{\left( \frac{x^{m}}{\omega^{m}}+\omega^{m}\right)^{\rho}}\right]\left(\frac{u}{ \omega}\right),\; u,\omega>0,\nonumber
\end{align}
and,
\begin{align}\label{M_L duality2}
\mathcal{M}_{\rho,m}\left[\left( \frac{x^{m}}{\omega^{m}}+\omega^{m}\right)^{\rho}e^{v\omega/x}f(x)\right](u,v,\omega)= \mathcal{L}\left[f(\omega x)\right](u),\; u,v,\omega>0.
\end{align}
Further, the Sumudu transform defined by (see, e.g., \cite{IT8,IT9,IT10})
\begin{align}\label{Sumudu IT}
\mathcal{S}[f(x)](\omega)=\int_{0}^{\infty}e^{-x} f(\omega x) dx,\; \omega>0
\end{align}
is a special case of the $\mathcal{M}$-transform \eqref{newIT}, where 
\[\mathcal{S}[f(x)](\omega)=\mathcal{M}_{0,m}[f(x)](1,0,\omega),\;\omega>0,\]
and, 
\begin{align*}
\mathcal{M}_{\rho, m}[f(x)](u,v,\omega)=
\frac{1}{u} \mathcal{S}\left[\frac{e^{-v\omega/x}f(x)}{\left(\frac{x^{m}}{\omega^{m}}+\omega^{m}\right)^{\rho}} \right] \left(\frac{\omega}{u}\right),\; u,v,\omega>0.
\end{align*}
Another special case of the integral transform \eqref{newIT}, when $u=v=0$, is a generalization of the Stieltjes transform, which was studied in (for example \cite{IT026,IT027}).
\begin{align}\label{Stieltjes IT}
\mathcal{S}t_{\rho}[f(x)](\omega)=\int_0^\infty \frac{f(x)}{\left(x +\omega\right)^\rho}dx,
\end{align}
with $\rho\in\mathbb C,\,\operatorname{Re}\rho>0$,  $\omega\in\mathbb C\setminus(-\infty,0)$.

Furthermore, important connections of the $\mathcal{M}$-transform with  well-known integral transforms will be discussed. These transforms include:\\
$\bullet$ The Borel-D\v zrbashjan transform defined (See, e.g., \cite{IT23}) as
\begin{align}\label{BorelIT}
\mathfrak{B}_{\nu,\mu}\left[f(x)\right](s)=\nu s^{\nu\mu-1} \int_0^\infty e^{-s^\nu x^{\nu}} x^{\nu\mu-1} f(x)dx,\,\, \nu,\mu>0.
\end{align}
$\bullet$ The $\bf{H}$-transform defined \cite{IT24} as
\begin{align}\label{HIT}
{\bf{H}}[f(x)](t)=\int_0^\infty H^{m,n}_{r,s}\left[xt\left\vert  \begin{array}{l}(a_i,\alpha_i)_{1,r}\\(b_j,\beta_j)_{1,s}\end{array} \right.\right]f(x)dx,
\end{align}
where $m,n,r,s$ are positive integers such that $0\leq m\leq s$ and $0\leq n\leq r$; $a_i, b_j \in\mathbb C$, $\alpha_i, \beta_j \in \mathbb R^+$ $(1\leq i\leq r; 1\leq j\leq s)$ and
\begin{align}\label{HFunction}
H^{m,n}_{r,s}\left[z\left\vert {\begin{array}{l}(a_i,\alpha_i)_{1,r}\\(b_j,\beta_j)_{1,s}\end{array}} \right.\right]
=\frac{1}{2\pi i}\int_\mathcal{C} \mathcal{H}^{m,n}_{r,s}(\theta) z^{-\theta}d\theta,
\end{align}
is the $H$-function defined in terms of a Mellin-Barnes type integral over a suitable contour $\mathcal C$, with
\begin{align}\label{HKernel}
\mathcal{H}^{m,n}_{r,s}(\theta)= \frac{\prod_{j=1}^{m}\Gamma(b_j +\beta_j \theta)}{\prod_{j=m+1}^{s}\Gamma(1-b_j -\beta_j \theta)} \frac{\prod_{i=1}^{n}\Gamma(1-a_j -\alpha_j \theta)}{\prod_{i=n+1}^{r}\Gamma(a_i +\alpha_i \theta)}.
\end{align}
Here, an empty product, if it exists, is taken to equal 1.\\
$\bullet$ The Mellin transform defined \cite{IT8,IT9} by
\begin{align}\label{Mellin-IT}
\mathfrak{M}[f(x)](z)=&\int_0^\infty x^{z-1}f(x)dx,\; \operatorname{Re} z>0.
\end{align}

In spirit of the above details, the $\mathcal M$-transform given by \eqref{newIT} seems to be worthy deserving to study because it has useful connections with the aforelisted integral transforms.

The present work is organized as follows. In Section 2, we discuss the existence of the $\mathcal{M}$-transform, further basic properties are studied. A new extension of the $H$-function is introduced, we call it an extended $H$-function. The relation of the $\mathcal{M}$-transform and this $H$-function is established. Also, the $\mathcal{M}$-transforms of differential derivatives are obtained. In Section 3, integral identities involving the $\mathcal{M}$-transform are obtained. Relationships of the $\mathcal{M}$-transform with the Borel-D\v{z}rbashjan transform, the Mellin transform and the $\bf{H}$-transform. Further, convolution and inversion formulas of the $\mathcal{M}$-transform are given. In Section 4, we present illustrative examples that can show the applicability of the $\mathcal{M}$-transform \eqref{newIT} in solving initial boundary value problems. Finally, we give concluding remarks in Section 5.

\section{Fundamental properties of the $\mathcal{M}$-transform}\label{sec2}

First, we study the existence of $\mathcal{M}$-transform defined by \eqref{newIT}. 
\begin{theorem}[Existence Theorem]\label{existenceTheorem}
	Let $f(x)$ be a continuous function ( or piecewise continuous) in $(0,\infty)$ satisfying
	\begin{align}\label{existencecodition}
	\lvert f(x)\rvert \leq K x^{m\operatorname{Re}\rho} e^{\frac{x}{\beta}}\; {\text{for all}} \;\; x>T,		
	\end{align}
	where $K,T$ and $\beta$ are positive constants. Then,  $\mathcal{M}_{\rho,m}[f(x)](u,v,\omega)$ given by \eqref{newIT} exists for all $u,v$ and $\omega$ such that 
	$$\omega\in(0,\mu),\; \operatorname{Re}u>\frac{\mu}{\beta},$$ 
	for some positive constant $\mu$. Furthermore, the integral \eqref{newIT} converges uniformly whenever $\operatorname{Re}u\geq a>\frac{\mu}{\beta}$ holds.
	\begin{proof}
		Using \eqref{existencecodition}, we get
		\begin{align}\label{existencecodition-1}
		\lvert f(\omega x)\rvert \leq K (\omega x)^{m\operatorname{Re}\rho} e^{\frac{\omega x}{\beta}},		
		\end{align}
		which, in view of
		\[\left\vert e^{-ux-v/x}\right\vert\leq e^{-x\operatorname{Re} u},\; x>0,\]
		leads to
		\begin{align*}
		\left\vert \mathcal{M}_{\rho,m}[f(x)](u,v,\omega)\right\vert & \leq \int_0^\infty \frac{e^{-x \operatorname{Re}u}}{\left(x^m +\omega^m \right)^{\operatorname{Re}\rho}}\lvert f(\omega x)\rvert dx\\
		%
		%
		&\leq K \omega^{m\operatorname{Re}\rho} \int_0^\infty e^{-\left(\operatorname{Re}u-\frac{\omega}{\beta}\right)x} x^{m\operatorname{Re}\rho}dx,
		%
		\end{align*}
		the last integral and hence the integral \eqref{newIT} exist if $\operatorname{Re}u>\frac{\mu}{\beta}$ and $0<\omega<\mu$.\\
That the integral \eqref{newIT} converges uniformly, follows directly from the well-known Weierstrass's test.
\end{proof}
\end{theorem}
Next, in view of Theorem \ref{existenceTheorem} we give the following basic properties of the $\mathcal{M}$-transform \eqref{newIT}.
\begin{theorem}[Scaling property]\label{scalingproperty} 
		Let $f$ satisfy the condition \eqref{existencecodition}. Then,
	\begin{align}\label{scalingpropertyEq}
	\mathcal{M}_{\rho,m}[f(\alpha^2 x)](u,v,\omega)=\alpha^{m\rho-1} \mathcal{M}_{\rho,m}[f(x)](\frac{u}{\alpha},\alpha v,\alpha\omega),\; \alpha>0,
	\end{align}
where $\rho\in\mathbb C,\,\operatorname{Re}\rho>0$, $m\in\mathbb N$ and $u,v\in\mathbb C$, and $\omega\in\mathbb R^+$.
	\begin{proof}
		This result follows directly from noting that
		\begin{align*}
		\begin{split}
		\mathcal{M}_{\rho,m}[f(\alpha^2 x)](u,v,\omega)=& \int_0^\infty \frac{e^{-\left(ux+\frac{v}{x}\right)}}{\left(x^m +\omega^m\right)^\rho}f(\alpha^2 \omega x)dx\\
		=&\alpha^{m\rho-1} \int_0^\infty \frac{e^{-\left(\frac{u}{\alpha}y+\frac{\alpha v}{y}\right)}}{\left((y)^m +(\alpha\omega)^m\right)^\rho}f(\alpha\omega y)dy,\quad (y=\alpha x)\\
		=& 	\alpha^{m\rho-1} \mathcal{M}_{\rho,m}[f(x)](\frac{u}{\alpha},\alpha v,\alpha\omega).
		\end{split}
		\end{align*}
	\end{proof}
\end{theorem}
Operting \eqref{newIT} on the function $\left(\frac{x^m}{\omega^m}+\omega^m\right)^\eta f(x)$, gives
	\[\mathcal{M}_{\rho,m}\big[\left(\frac{x^m}{\omega^m}+\omega^m\right)^\eta f(x)\big](u,v,\omega)=\int_0^\infty \frac{e^{-\left(ux+\frac{v}{x}\right)}}{\left(x^m +\omega^m\right)^{\rho-\eta}}f(\omega x)dx,	\]
which leads directly to the following results.
\begin{theorem}[Elimination property]\label{Elimination property}
	If \(\lvert f(x)\rvert \leq K x^{m\operatorname{Re}(\rho-\eta)}e^{\frac{x}{\beta}}\) with \(\operatorname{Re}(\rho-\eta)\geq 0,\,K,\beta>0\), then
	\begin{align}\label{Elimination propertyEq}
	\mathcal{M}_{\rho,m}\big[\left(\frac{x^m}{\omega^m}+\omega^m\right)^\eta f(x)\big](u,v,\omega)= \mathcal{M}_{\rho-\eta,m}[f(x)](u,v,\omega).
	\end{align}
	For \(\eta=\rho\), \eqref{Elimination propertyEq} reduces to
	\begin{align}\label{Elimination propertyEq2}
	\mathcal{M}_{\rho,m}\big[\left(\frac{x^m}{\omega^m}+\omega^m\right)^\rho f(x)\big](u,v,\omega)= \mathcal{M}_{0,m}[f(x)](u,v,\omega).
	\end{align}
\end{theorem}
\begin{remark}
	The \(\mathcal{M}\)-transform \eqref{newIT} can be used to eliminate the factor \(\left(\frac{x^m}{\omega^m}+\omega^m\right)^\rho\) and reduce it to a simple form, the natural transform defined by \eqref{naturalIT}.
\end{remark}
For further investigations, we need to recall the extended gamma function defined (See, e.g., \cite{IT25}) as
\begin{align}\label{extended gamma}
\Gamma_b (z)=\int_0^\infty t^{z-1} e^{-t-b/t}dt, \,\,\operatorname{Re} z>0,\,\operatorname{Re} b>0.
\end{align}
Simple computations show that
\begin{align}\label{EulerIntegral-extend}
\int_0^\infty x^{z-1} e^{-\sigma x-\frac{b}{x}}dx=\frac{\Gamma_{\sigma b}(z)}{\sigma^z},\quad \operatorname{Re} z>0,\, \operatorname{Re}\sigma>0,\, \operatorname{Re}b>0.
\end{align}
%
In particular, when $b=0$ we get the Euler integral
\begin{align}\label{EulerIntegral}
\int_0^\infty x^{z-1} e^{-\sigma x}dx=\frac{\Gamma(z)}{\sigma^z},\quad \operatorname{Re} z>0,\, \operatorname{Re}\sigma>0.
\end{align}
Also, in spirit of \eqref{HFunction}, \eqref{HKernel} and \eqref{extended gamma} we introduce the following extended $H$-function
\begin{align}\label{extend HFunction}
H^{2,1}_{1,2}(z;b)&=H^{2,1}_{1,2}\left[z;b\left\vert {\begin{array}{l}(a,\alpha)\\(b_1,\beta_1)_{b},(b_2,\beta_2) \end{array}}\right.\right]\nonumber\\
&=\frac{1}{2\pi i}\int_{c-i\infty}^{c+i\infty} \Gamma(1-a-\alpha t) \Gamma_b (b_1 +\beta_1 t) \Gamma(b_2 +\beta t)z^{-t}dt,\; c\geq 0,
\end{align}
where $\Gamma_b$ is the extend gamma function \eqref{extended gamma} and $\Gamma$ is the classical gamma function.
\begin{align*}
\Gamma(z)=\int_0^\infty t^{z-1} e^{-t}dt, \,\,\operatorname{Re} z>0.
\end{align*}
Since, for $m\in\mathbb N$,
\begin{align*}
\int_{c-i\infty}^{c+i\infty}& \Gamma(1-a-\alpha t) \Gamma_b (b_1 +\beta_1 t) \Gamma(b_2 +\beta t)z^{-mt}dt\\
&=\frac{1}{m}\int_{c/m -i\infty}^{c/m +i\infty} \Gamma(1-a-\frac{\alpha}{m} t) \Gamma_b (b_1 +\frac{\beta_1}{m} t) \Gamma(b_2 +\frac{\beta}{m} t)z^{-t}dt,
\end{align*}
then, the extended $H$-function \eqref{extend HFunction} satisfies the following identity
%
\begin{align}\label{lem1/m}
H^{2,1}_{1,2}\left[z^m;b\left\vert {\begin{array}{l}(a,\alpha)\\(b_1,\beta_1)_{b},(b_2,\beta_2) \end{array}} \right.\right]
=\frac{1}{m}
H^{2,1}_{1,2}\left[z;b\left\vert {\begin{array}{l}(a,\frac{\alpha}{m})\\(b_1,\frac{\beta_1}{m})_{b},(b_2,\frac{\beta_2}{m}) \end{array}}\right.\right].
\end{align}
%
%
The following theorem provides the \(\mathcal{M}\)-images of power and exponential functions under the \(\mathcal{M}\)-transform \eqref{newIT}. 
\begin{theorem}\label{M_PowerExponential}
	For any \(\rho\in\mathbb C\), with \(\operatorname{Re} \rho>0\), \(\lambda>0\) and \(a\geq 0\) we have
	\begin{align}\label{M-powerFunction}
	\mathcal{M}_{\rho,m}\left[x^{\lambda-1}\right](u,v,\omega)= \frac{\omega^{\lambda-m\rho-1}u^{-\lambda}}{m\Gamma(\rho)} H^{2,1}_{1,2} \left[u\omega\left\vert{\begin{array}{l}
		\left(1,\frac{1}{m}\right)\\ \left(\lambda,1\right)_{uv},\left(\rho,\frac{1}{m}\right)
		\end{array}}\right.\right],
	\end{align}
	\begin{align}\label{M-exponentialFunction}
	\mathcal{M}_{\rho.m}\left[e^{-a x}\right](u,v,\omega)= \frac{\omega^{-m\rho}}{m(u+a\omega)\Gamma(\rho)} H^{2,1}_{1,2} \left[\omega(u+a\omega)\left\vert{\begin{array}{l}
		\left(1,\frac{1}{m}\right)\\ \left(1,1\right)_{(u+a\omega)v},\left(\rho,\frac{1}{m}\right)
		\end{array}}\right.\right],
	\end{align}
	and
\begin{align}\label{M-p-expFun}
\mathcal{M}_{\rho.m}\left[x^{\lambda-1} e^{-a x}\right](u,v,\omega)= \frac{\omega^{\lambda-m\rho-1}}{m(u+a\omega)^{\lambda}\Gamma(\rho)} H^{2,1}_{1,2} \left[\omega(u+a\omega)\left\vert{\begin{array}{l}
	\left(1,\frac{1}{m}\right)\\ \left(\lambda,1\right)_{uv},\left(\rho,\frac{1}{m}\right)
	\end{array}}\right.\right],
\end{align}
where $H^{2,1}_{1,2}$ is the extended $H$-function defined by \eqref{extend HFunction}.
\end{theorem}
\begin{proof}
	$\bullet$ For the first conclusion \eqref{M-powerFunction}, it is clear that
	\begin{align}\label{M-powerFunction2}
	\mathcal{M}_{\rho,m}\left[x^{\lambda-1}\right](u,v,\omega)=& \int_0^\infty \frac{e^{-ux-\frac{v}{x}}}{\left(x^m +\omega^m\right)^\rho}(\omega x)^{\lambda-1}dx\nonumber\\
	=&\frac{\omega^{\lambda-1}}{\Gamma(\rho)}\int_0^\infty e^{-ux-\frac{v}{x}}x^{\lambda-1} \left[\int_0^\infty s^{\rho-1} e^{-\left(x^m +\omega^m\right)s}ds\right]dx\nonumber\\
	=&\frac{\omega^{\lambda-1}}{\Gamma(\rho)}\int_0^\infty s^{\rho-1} e^{-\omega^m s} \left[\int_0^\infty e^{-ux-\frac{v}{x}-x^m s}x^{\lambda-1}dx\right]ds.
	\end{align}
	Now, in order to evaluate the integral
	\begin{align}\label{I(s)1}
	I(s)=\int_0^\infty e^{-ux-\frac{v}{x}-x^m s}x^{\lambda-1}dx,\, s>0
	\end{align}
	we apply the Mellin transform, defined by \eqref{Mellin-IT}.
	For \(\operatorname{Re} z>0,\,\operatorname{Re}(\lambda-mz)>0\), and \(\operatorname{Re} uv>0\), we have
	\begin{align}\label{MellinI}
	\mathfrak{M}[I(s)](z)=&\int_0^\infty x^{\lambda-1} e^{-ux-\frac{v}{x}}\left[\int_0^\infty  s^{z-1} e^{-x^m s}ds\right]dx\nonumber\\
	=&\Gamma(z) \int_0^\infty x^{\lambda-mz-1} e^{-ux-\frac{v}{x}}dx\nonumber\\
	=&\frac{\Gamma(z) \Gamma_{uv}(\lambda-mz)}{u^{\lambda-mz}},
	\end{align}
	here, \eqref{EulerIntegral-extend} is used. Thus, the integral $I(s)$ defined by \eqref{I(s)1}, can be evaluated by applying the inverse Mellin transform to \eqref{MellinI}, then we get
	\begin{align}\label{I(s)2}
	I(s)=\frac{u^{-\lambda}}{2\pi i}\int_{-i\infty}^{i\infty} \Gamma(z) \Gamma_{uv}(\lambda-mz) u^{mz}s^{-z}dz,\, s>0
	\end{align}
	Substituting \eqref{I(s)2} into \eqref{M-powerFunction2}, and interchanging the order of integration, give
	\begin{align}\label{M-powerFunction3}
	\mathcal{M}_{\rho,m}\left[x^{\lambda-1}\right](u,v,\omega)=& \frac{\omega^{\lambda-1}u^{-\lambda}}{2\pi i\Gamma(\rho)}\int_0^\infty s^{\rho-1} e^{-\omega^m s} \int_{-i\infty}^{i\infty} \Gamma(z) \Gamma_{uv}(\lambda-mz)u^{mz} s^{-z}dz ds\nonumber\\
	=& \frac{\omega^{\lambda-1}u^{-\lambda}}{2\pi i\Gamma(\rho)} \int_{-i\infty}^{i\infty} \Gamma(z) \Gamma_{uv}(\lambda-mz)u^{mz} \int_0^\infty s^{\rho-z-1} e^{-\omega^m s}ds dz\nonumber\\
	=& \frac{\omega^{\lambda-m\rho-1}u^{-\lambda}}{2\pi i\Gamma(\rho)} \int_{-i\infty}^{i\infty} \Gamma(z) \Gamma_{uv}(\lambda-mz) \Gamma(\rho-z) \left(\omega^m u^m\right)^z dz\nonumber\\
	=& \frac{\omega^{\lambda-m\rho-1}u^{-\lambda}}{2\pi i\Gamma(\rho)} \int_{-i\infty}^{i\infty} \Gamma(-z) \Gamma_{uv}(\lambda+mz) \Gamma(\rho+z) \left(\omega^m u^m\right)^{-z} dz\nonumber\\
	=& \frac{\omega^{\lambda-m\rho-1}u^{-\lambda}}{\Gamma(\rho)} H^{2,1}_{1,2} \left[\omega^m u^m \left\vert {\begin{array}{l}
		\left(1,1\right)\\ \left(\lambda,m\right)_{uv},\left(\rho,1\right)
		\end{array}}\right.\right]\nonumber\\
	=& \frac{\omega^{\lambda-m\rho-1}u^{-\lambda}}{m\Gamma(\rho)} H^{2,1}_{1,2} \left[\omega u \left\vert {\begin{array}{l}
		\left(1,\frac{1}{m}\right)\\ \left(\lambda,1\right)_{uv},\left(\rho,\frac{1}{m}\right)
		\end{array}}\right.\right],
	\end{align} 
	the last equality follows directely from \eqref{lem1/m}. This proves the required conclusion \eqref{M-powerFunction}.
	
	$\bullet$ For the second conclusion \eqref{M-exponentialFunction}, we have
	\begin{align}\label{M-expFun2}
	\mathcal{M}_{\rho.m}\left[e^{-a x}\right](u,v,\omega)=&\int_0^\infty \frac{e^{-ux-\frac{v}{x}}}{\left(x^m +\omega^m\right)^\rho} e^{-a\omega x} dx\nonumber\\
	=&\frac{1}{\Gamma(\rho)}\int_0^\infty e^{-ux-\frac{v}{x}} e^{-a\omega x} \int_0^\infty s^{\rho-1} e^{-\left(x^m +\omega^m\right)s}ds dx\nonumber\\
	=&\frac{1}{\Gamma(\rho)}\int_0^\infty s^{\rho-1} e^{-\omega^m s} \int_0^\infty e^{-ux-\frac{v}{x}-a\omega x-x^m s}dx ds.
	\end{align}
	Once again, we apply the Mellin transform to evaluate the following integral
	\begin{align}\label{J(s)1}
	J(s)=\int_0^\infty e^{-(u+a\omega)x-\frac{v}{x}-x^m s}dx.
	\end{align}
	Doing so, gives
	\begin{align}
	\mathfrak{M}[J(s)](z)=&\int_0^\infty s^{z-1}\int_0^\infty e^{-(u+a\omega)x-\frac{v}{x}-x^m s}dx ds\nonumber\\
	%
	%
	%
	=& \frac{\Gamma(z)}{(u+a\omega)^{1-mz}}\Gamma_{(u+a\omega)v}(1-mz).
	\end{align}
	Applying the inverse Mellin transform, gives
	\begin{align}\label{J(s)2}
	J(s)=\frac{(u+a\omega)^{-1}}{2\pi i} \int_{-i\infty}^{+i\infty} \Gamma(z) \Gamma_{(u+a\omega)v}(1-mz) (u+a\omega)^{mz} s^{-z}dz
	\end{align}
	Substituting \eqref{J(s)2} into \eqref{M-expFun2}, leads to
	\begin{align}
	\mathcal{M}_{\rho,m}\left[e^{-ax}\right](u,v,\omega)=&\frac{(u+a\omega)^{-1}}{2\pi i\Gamma(\rho)} \int_0^\infty s^{\rho-1} e^{-\omega^m s}\nonumber\\
	&\times \int_{-i\infty}^{+i\infty} \Gamma(z) \Gamma_{(u+a\omega)v}(1-mz) (u+a\omega)^{mz} s^{-z}dz ds\nonumber\\
	=&\frac{(u+a\omega)^{-1}}{2\pi i\Gamma(\rho)} \int_{-i\infty}^{+i\infty} \Gamma(z) \Gamma_{(u+a\omega)v}(1-mz) (u+a\omega)^{mz}\nonumber\\
	&\times  \int_0^\infty s^{\rho-z-1} e^{-\omega^m s}ds dz\nonumber\\
	=&\frac{(u+a\omega)^{-1}\omega^{-m\rho}}{2\pi i\Gamma(\rho)} \int_{-i\infty}^{+i\infty} \Gamma(z) \Gamma_{(u+a\omega)v}(1-mz)\nonumber\\
	&\times\Gamma(\rho-z) \left[(u+a\omega)\omega\right]^{mz} dz\nonumber\\
	=&\frac{\omega^{-m\rho}}{(u+a\omega)\Gamma(\rho)} H^{2,1}_{1,2} \left[\omega^m (u+a\omega)^m \left\vert {\begin{array}{l}
		\left(1,1\right)\\ \left(1,m\right)_{(u+a\omega)v},\left(\rho,1\right)
		\end{array}}\right.\right],\nonumber
	\end{align}
	which in view of \eqref{lem1/m} reduces directly to the conclusion \eqref{M-exponentialFunction}.

	$\bullet$ The third conclusion \eqref{M-p-expFun} of Theorem \ref{M_PowerExponential} can be obtained directly from \eqref{M-powerFunction} and noting that 
	\begin{align}\label{M-p-expFun2}
	\mathcal{M}_{\rho,m}\left[x^{\lambda-1} e^{-a x}\right](u,v,\omega)=
	\mathcal{M}_{\rho,m}\left[x^{\lambda-1}\right](u+a\omega,v,\omega).
	\end{align}
\end{proof}
\begin{remark}
	If we set $\lambda=1$ in \eqref{M-powerFunction}, or $a=0$ in \eqref{M-exponentialFunction}, we get
	\begin{align}\label{M_1}
	\mathcal{M}_{\rho.m}\left[1\right](u,v,\omega)= \frac{\omega^{\lambda-m\rho-1}}{m u\Gamma(\rho)} H^{2,1}_{1,2} \left[u\omega \left\vert {\begin{array}{l}
		\left(1,\frac{1}{m}\right)\\ \left(1,1\right)_{uv},\left(\rho,\frac{1}{m} \right)
		\end{array}}\right.\right].
	\end{align}
\end{remark}	
Further results, that can be directly computed are summarized in Table 1. 
\begin{table}[h!]\label{table1}
	\begin{center}
		\caption{Further $\mathcal{M}_{\rho,m}$-images}
		\label{tab:table1}
		\begin{tabular}{c|c} 
		\hline\hline
		\textbf{Function} & \textbf{New integral transform}\\
				$f(x)$ & $\mathcal{M}_{\rho,m}[f(x)](u,v,\omega)$ \\
			\hline\\
			$x^n f(x)$ & $(-1)^n \omega^n \frac{\partial^n}{\partial u^n}\mathcal{M}_{\rho,m}[f(x)](u,v,\omega)$ \\ \\
			$x^n f(x)$ & $\omega^{n} \int_v^\infty \int_{s_n}^\infty \cdots\int_{s_2}^\infty \mathcal{M}_{\rho,m}[f(x)](u,s_1,\omega)ds_1 \cdots ds_{n-1} ds_n$ \\ \\
			$\frac{f(x)}{x^n}$ & $\omega^{-n} \int_u^\infty \int_{s_n}^\infty \cdots\int_{s_2}^\infty \mathcal{M}_{\rho,m}[f(x)](s_1,v,\omega)ds_1 \cdots ds_{n-1} ds_n$ \\ \\
			$\frac{f(x)}{x^n}$ & $(-1)^n\omega^{-n} \frac{\partial^n}{\partial v^n} \mathcal{M}_{\rho,m}[f(x)](u,v,\omega)$\\ \\
			$e^{-a/x}f(x)$ & $\mathcal{M}_{\rho,m}[f(x)](u,v+a/\omega ,\omega)$\\ \\
			$e^{-ax}f(x)$ & $\mathcal{M}_{\rho,m}[f(x)](u+a\omega,v,\omega)$\\
			\hline\hline
		\end{tabular}
	\end{center}
\end{table}
\newpage
Next, we give the $\mathcal{M}$-transform of derivatives.
\begin{theorem}[$\mathcal{M}$-transform of derivatives]\label{M_derivatives}
	Let the assumptions of Theorem \ref{existenceTheorem} be satisfied for the $n^{th}$ derivative of a function $f(x)$. Then,
	\begin{align}\label{M-fn}
	\mathcal{M}_{\rho,m}\left[f^{(n)}(x)\right](u,v,\omega)=& \frac{u^n}{\omega^n} \mathcal{M}_{\rho,m}\left[f(x)\right](u,v,\omega) - \delta_{v,0} \sum_{k=0}^{n-1} \frac{u^k}{\omega^{m\rho+k+1}} f^{(n-k-1)}(0)\nonumber\\
	&\, -v\omega \sum_{k=0}^{n-1} \frac{u^k}{\omega^{k}} \mathcal{M}_{\rho,m}\left[x^{-2}f^{(n-k-1)}(x)\right](u,v,\omega) \nonumber\\
	&\,+\frac{m\rho}{\omega^m} \sum_{k=0}^{n-1} \frac{u^k}{\omega^{k}} \mathcal{M}_{\rho+1,m}\left[x^{m-1}f^{(n-k-1)}(x)\right](u,v,\omega),
	\end{align}
	where $\delta_{v,0}$ is the delta function, defined as
	\[\delta_{v,0}=\begin{cases}1,\,& v=0,\\ 0, \, &v\neq 0.\end{cases}\]
\end{theorem}
\begin{proof}
	
	$\bullet$ For $n=1$, if $f(x)$ and $f'(x)$ satisfy the conditions of Theorem \ref{existenceTheorem}, then
	\begin{align*}
	\mathcal{M}_{\rho,m}\left[f'(x)\right](u,v,\omega)=& \int_0^\infty \frac{e^{-ux-\frac{v}{x}}}{\left(x^m +\omega^m\right)^\rho} f'(\omega x)dx=\frac{1}{\omega} \int_0^\infty \frac{e^{-ux-\frac{v}{x}}}{\left(x^m +\omega^m\right)^\rho} df(\omega x)\\
	=&\frac{1}{\omega} \left[\left.\frac{e^{-ux-\frac{v}{x}}f(\omega x)}{\left(x^m +\omega^m\right)^\rho}\right\vert_0^\infty -\int_0^\infty f(\omega x) \frac{\partial}{\partial x} \frac{e^{-ux-\frac{v}{x}}}{\left(x^m +\omega^m\right)^\rho} dx
	\right].
	\end{align*}
	In view of the asymptotic property \eqref{existencecodition}, for \(\operatorname{Re} u>\frac{\omega}{\beta}\) and \(\operatorname{Re} v>0\) when \(v\in\mathbb C\setminus \{0\}\), we have
	\begin{align*}
	\left\vert \frac{e^{-ux-\frac{v}{x}}f(\omega x)}{\left(x^m +\omega^m\right)^\rho} \right\vert \leq K \omega^{m\operatorname{Re}\rho} e^{-\left(\operatorname{Re} u -\omega/\beta\right)x},
	\end{align*}
	which shows that
	\begin{align*}
	\lim_{x\to\infty} \frac{e^{-ux-\frac{v}{x}}f(\omega x)}{\left(x^m +\omega^m\right)^\rho}=0,\,\, \lim_{x\to 0} \frac{e^{-ux-\frac{v}{x}}f(\omega x)}{\left(x^m +\omega^m\right)^\rho}=\frac{f(0)}{\omega^{m\rho}}\delta_{v,0}.
	\end{align*}
	On the other hand, we have
	\begin{align*}
	\frac{\partial}{\partial x} \frac{e^{-ux-\frac{v}{x}}}{\left(x^m +\omega^m\right)^\rho}= -\left(u-\frac{v}{x^2}\right) \frac{e^{-ux-\frac{v}{x}}}{\left(x^m +\omega^m\right)^\rho} -m\rho x^{m-1}\frac{e^{-ux-\frac{v}{x}}}{\left(x^m +\omega^m\right)^{\rho+1}}.
	\end{align*}
	Combining these results together, gives
	\begin{align}\label{M_F'}
	\mathcal{M}_{\rho,m}\left[f'(x)\right](u,v,\omega)=&-\frac{f(0)}{\omega^{m\rho+1}}\delta_{v,0}+ \frac{u}{\omega} \int_0^\infty \frac{e^{-ux-\frac{v}{x}}}{\left(x^m +\omega^m\right)^\rho} f(\omega x)dx \nonumber\\
	&\, -\frac{v}{\omega}\int_0^\infty \frac{e^{-ux-\frac{v}{x}}}{\left(x^m +\omega^m\right)^\rho} \frac{f(\omega x)}{x^2}dx\nonumber\\
	&\, +\frac{m\rho}{\omega}\int_0^\infty \frac{e^{-ux-\frac{v}{x}}}{\left(x^m +\omega^m\right)^\rho+1}x^{m-1} f(\omega x) dx\nonumber\\
	=& \frac{u}{\omega} \mathcal{M}_{\rho,m}\left[f(x)\right](u,v,\omega)-\frac{f(0)}{\omega^{m\rho+1}}\delta_{v,0}\nonumber\\
	&\, -v\omega \mathcal{M}_{\rho,m}\left[x^{-2}f(x)\right](u,v,\omega)\nonumber\\
	&\, + \frac{m\rho}{\omega^m} \mathcal{M}_{\rho+1,m}\left[x^{m-1}f(x)\right](u,v,\omega).
	\end{align}
	$\bullet$ For the induction, assume that \eqref{M-fn} holds for some integer $n$.  
	Then, for $n+1$, on the basis of \eqref{M_F'}, one gets
	\begin{align}\label{M_Fn+1}
	\mathcal{M}_{\rho,m}\left[f^{(n +1)}(x)\right](u,v,\omega)
	=& \frac{u}{\omega} \mathcal{M}_{\rho,m}\left[f^{(n)}(x)\right](u,v,\omega)-\frac{f^{(n)}(0)}{\omega^{m\rho+1}}\delta_{v,0}\nonumber\\
	&\, -v\omega \mathcal{M}_{\rho,m}\left[x^{-2}f^{(n)}(x)\right](u,v,\omega)\nonumber\\
	&\, + \frac{m\rho}{\omega^m} \mathcal{M}_{\rho+1,m}\left[x^{m-1}f^{(n)}(x)\right](u,v,\omega).
	\end{align}
	Plugging \eqref{M-fn} into \eqref{M_Fn+1}, leads to
	\begin{align*}
	\mathcal{M}_{\rho,m}\left[f^{(n +1)}(x)\right](u,v,\omega)
	=& \frac{u^{n+1}}{\omega^{n+1}} \mathcal{M}_{\rho,m}\left[f(x)\right](u,v,\omega) \nonumber\\
	&\, - \delta_{v,0} \sum_{k=0}^{n-1} \frac{u^{k+1}}{\omega^{m\rho+k+2}} f^{(n -k-1)}(0)\nonumber\\
	&\, -v\omega \sum_{k=0}^{n -1} \frac{u^{k+1}}{\omega^{k+2}} \mathcal{M}_{\rho,m}\left[x^{-2}f^{(n -k-1)}(x)\right](u,v,\omega) \nonumber\\
	&\,+\frac{m\rho}{\omega^m} \sum_{k=0}^{n-1} \frac{u^{k+1}}{\omega^{k+2}} \mathcal{M}_{\rho+1,m}\left[x^{m-1}f^{(n -k-1)}(x)\right](u,v,\omega)\nonumber\\
	&\, -\frac{f^{(n)}(0)}{\omega^{m\rho+1}}\delta_{v,0} -v\omega \mathcal{M}_{\rho,m}\left[x^{-2}f^{(n)}(x)\right](u,v,\omega)\nonumber\\
	&\, + \frac{m\rho}{\omega^m} \mathcal{M}_{\rho+1,m}\left[x^{m-1}f^{(n)}(x)\right](u,v,\omega),
	\end{align*}
	which can be rewritten as
	\begin{align*}
	\mathcal{M}_{\rho,m}\left[f^{(n +1)}(x)\right](u,v,\omega)
	=& \frac{u^{n+1}}{\omega^{n+1}} \mathcal{M}_{\rho,m}\left[f(x)\right](u,v,\omega) \nonumber\\
	&\, - \delta_{v,0} \sum_{k=0}^{n} \frac{u^{k}}{\omega^{m\rho+k+1}} f^{(n -k)}(0)\nonumber\\
	&\, -v\omega \sum_{k=0}^{n} \frac{u^{k}}{\omega^{k+1}} \mathcal{M}_{\rho,m}\left[x^{-2}f^{(n -k)}(x)\right](u,v,\omega) \nonumber\\
	&\,+\frac{m\rho}{\omega^m} \sum_{k=0}^{n} \frac{u^{k}}{\omega^{k+1}} \mathcal{M}_{\rho+1,m}\left[x^{m-1}f^{(n -k)}(x)\right](u,v,\omega).
	\end{align*}
	This completes the proof of Theorem \ref{M_derivatives}.
\end{proof}
\begin{remark}
	Many special cases can be deduced from \eqref{M-fn} as following
	\begin{itemize}
		\item[\textbf{1.}]  If we set $v=0$ in \eqref{M-fn}, we get (see, \cite{IT26})
		\begin{align*}
		\mathbb M_{\rho,m}\left[f^{(n)}(x)\right](u,\omega)=& \frac{u^n}{\omega^n} \mathbb M_{\rho,m}\left[f(x)\right](u,\omega) - \sum_{k=0}^{n-1} \frac{u^k}{\omega^{m\rho+k+1}} f^{(n-k-1)}(0)\nonumber\\
		&\,+\frac{m\rho}{\omega^m} \sum_{k=0}^{n-1} \frac{u^k}{\omega^{k}} \mathbb M_{\rho+1,m}\left[x^{m-1}f^{(n-k-1)}(x)\right](u,\omega),
		\end{align*}
		where $\mathbb M_{\rho,m}$ is the Srivastava-Luo-Raina integral transform defined by \eqref{SriIT}, where $v$ replaced by $\omega$.
		\item[\textbf{2.}] Setting $\rho=0, v=0$ in \eqref{M-fn}, gives
		\begin{align*}
		\mathcal{N}\left[f^{(n)}(x)\right](u,\omega)=& \frac{u^n}{\omega^n} \mathcal{N}\left[f(x)\right](u,\omega) - \sum_{k=0}^{n-1} \frac{u^k}{\omega^{k+1}} f^{(n-k-1)}(0),
		\end{align*}
		which can be rewritten as (see \cite{0IT20})
		\begin{align*}
		\mathcal{N}\left[f^{(n)}(x)\right](u,\omega)=& \frac{u^n}{\omega^n} \mathcal{N}\left[f(x)\right](u,\omega) - \sum_{k=0}^{n-1} \frac{u^{n-k-1}}{\omega^{n-k}} f^{(k)}(0),
		\end{align*}
		where $\mathcal{N}$ is the natural transform defined by \eqref{naturalIT}.
		\item[\textbf{3.}] If $\rho=v=0,\, u=1$, \eqref{M-fn} reduces to \cite{0IT16}
		\begin{align*}
		\mathcal{S}\left[f^{(n)}(x)\right](\omega)=& \frac{1}{\omega^n} \mathcal{S}\left[f(x)\right](\omega) - \sum_{k=0}^{n-1} \frac{1}{\omega^{n-k}} f^{(k)}(0),
		\end{align*}
		where $\mathcal{S}$ is the Sumudu transform defined by \eqref{Sumudu IT} .
		\item[\textbf{4.}] If $\rho=v=0,\, \omega=1$, \eqref{M-fn} reduces to \cite{IT8,IT9}
		\begin{align*}
		\mathcal{L}\left[f^{(n)}(x)\right](u)=& u^n \mathcal{L}\left[f(x)\right](u) - \sum_{k=0}^{n-1} u^{n-k-1} f^{(k)}(0),
		\end{align*}
		where $\mathcal{L}$ is the Laplace transform defined by \eqref{Laplace_T}.
	\end{itemize}
\end{remark}
Using Theorem \ref{Elimination property} the following consequences of the Theorem \ref{M_derivatives} can be achieved.
\begin{corollary}\label{corollary1}
	If $f^{(n)}$, the $n^{th}$-derivative of a function $f$, satisfies the assumption of Theorem \ref{existenceTheorem}, then for $\nu\in\mathbb C$ with $\operatorname{Re}\nu>0$ we have
	\begin{align}\label{M-nth derivative}
	\mathcal{M}_{\rho,m}\Big[\left( \frac{x^m}{\omega^m}+\omega^m\right)^\nu &f^{(n)}(x)\Big](u,v,\omega)= \frac{u^n}{\omega^n} \mathcal{M}_{\rho-\nu,m}\left[f(x)\right](u,v,\omega) \nonumber\\
	&\,- \delta_{v,0} \sum_{k=0}^{n-1} \frac{u^k}{\omega^{m(\rho-\nu)+k+1}} f^{(n-k-1)}(0)\nonumber\\
	&\, -v\omega \sum_{k=0}^{n-1} \frac{u^k}{\omega^{k}} \mathcal{M}_{\rho-\nu,m}\left[x^{-2}f^{(n-k-1)}(x)\right](u,v,\omega) \nonumber\\
	&\,+\frac{m(\rho-\nu)}{\omega^m} \sum_{k=0}^{n-1} \frac{u^k}{\omega^{k}} \mathcal{M}_{\rho-\nu+1,m}\left[x^{m-1}f^{(n-k-1)}(x)\right](u,v,\omega).
	\end{align}
\end{corollary}
\section{Integrals involving the $\mathcal{M}$-transform}
In this section, we present many integrals comprising the $\mathcal{M}$-transform \eqref{newIT}, well-known integral transforms and special functions. The identities established here include Parseval-type indentity, Mellin integral, $\bf{H}$-transform and Borel-D\v{z}rbashjan transform of the \(\mathcal{M}\)-transform. Also, convolution and inversion formulas are given.
\begin{theorem}[Parseval-type identity]\label{parseval}
	
	Let $f$ and $g$ satisfy the condition \eqref{existencecodition}. Then,
	\begin{align}\label{parsevalIdentity}
	\int_0^\infty \frac{f(\omega\tau) e^{-v/\tau}}{\left(\tau^m +\omega^m\right)^{\rho_1}}& \mathcal{M}_{\rho_2 ,m}[g(x)](\tau,v,\omega)d\tau \nonumber\\
	&= 
	\int_0^\infty \frac{g(\omega x) e^{-v/x}}{\left(x^m +\omega^m\right)^{\rho_2}} \mathcal{M}_{\rho_1 ,m} [f(\tau)](x,v,\omega)dx.
	\end{align}
	where $\rho_i\in\mathbb C,\,\operatorname{Re}\rho_i>0,\,(i=1,2)$, $m\in\mathbb N$ and $u,v\in\mathbb C$, and $\omega\in\mathbb R^+$
\end{theorem}
\begin{proof}
	Using \eqref{newIT} and applying the Fubini's theorem, give
	\begin{align*}
	\int_0^\infty& \frac{f(\omega\tau) e^{-v/\tau}}{\left(\tau^m +\omega^m\right)^{\rho_1}} \mathcal{M}_{\rho_2 ,m}[g(x)](\tau,v,\omega)d\tau\\
	&=\int_0^\infty \frac{f(\omega\tau) e^{-v/\tau}}{\left(\tau^m +\omega^m\right)^{\rho_1}} \int_0^\infty \frac{e^{-\tau x-v/x}}{\left(x^m +\omega^m\right)^{\rho_2}} g(\omega x) dx d\tau\\
	&=\int_0^\infty \frac{g(\omega x) e^{-v/x}}{\left(x^m +\omega^m\right)^{\rho_2}} \int_0^\infty \frac{e^{-\tau x-v/\tau}}{\left(\tau^m +\omega^m\right)^{\rho_1}} f(\omega\tau) d\tau dx\\
	&=\int_0^\infty \frac{g(\omega x) e^{-v/x}}{\left(x^m +\omega^m\right)^{\rho_2}} \mathcal{M}_{\rho_1 ,m} [f(\tau)](x,v,\omega)dx. 
	\end{align*}
	Which completes the proof of Theorem \ref{parseval}.
\end{proof}
\begin{remark}
	Setting $\rho_1 =\rho_2 =v=0$ in \eqref{parsevalIdentity}, leads to the Parseval-type identity for the natural transform \eqref{naturalIT} (See \cite{IT9,IT20})
	\begin{align}\label{parsevalNatural}
	\int_0^\infty f(\omega\tau) \mathcal{N}[g(x)](\tau,\omega)d\tau = 
	\int_0^\infty g(\omega x) \mathcal{N} [f(\tau)](x,\omega)dx.
	\end{align}
	Moreover, if we set $\omega=1$ in \eqref{parsevalNatural}, then we get the Parseval-type identity for the Laplace transform
	\begin{align*}
	\int_0^\infty f(\tau) \mathcal{L}[g(x)](\tau)d\tau = 
	\int_0^\infty g(x) \mathcal{L} [f(\tau)](x)dx.
	\end{align*}
\end{remark}
Applying the argument used in proving Theorem \ref{parseval}, we can prove the following result.
\begin{theorem}
	
	Let $f$ and $g$ satisfy the condition \eqref{existencecodition}. Then,
	\begin{align}
	\int_0^\infty \frac{f(\omega\tau) e^{-u\tau}}{\left(\tau^m +\omega^m\right)^{\rho_1}}& \mathcal{M}_{\rho_2 ,m}[g(x)](u,\tau,\omega)d\tau \nonumber\\
	& =
	\int_0^\infty \frac{g(\omega x) e^{-ux}}{\left(x^m +\omega^m\right)^{\rho_2}} \mathcal{M}_{\rho_1 ,m} [f(\tau)](u+\frac{1}{x},0,\omega)dx.
	\end{align}
\end{theorem}
Next, we obtain identities providing the relationships between the $\mathcal{M}$-transform \eqref{newIT} and some well-known integral transforms.
\begin{theorem}
	
	Let $f$ and $g$ satisfy the condition \eqref{existencecodition}. Then,
	\begin{align}\label{M-N}
	\int_0^\infty f(\omega u) \mathcal{M}_{\rho ,m}[e^{-ax}g(x)](u\omega,v,\omega)du = 
	\mathcal{M}_{\rho,m}\big[g(x)\mathcal{N}[f(u)](x,\omega)\big]
	(a\omega,v,\omega),\, a>0,
	\end{align}
	where $\mathcal{N}$ is the natural transform defined by \eqref{naturalIT}.
\end{theorem}
\begin{proof}
	Plugging \eqref{newIT} into the left-hand side of \eqref{M-N}, gives
	\begin{align*}
	\int_0^\infty f(\omega u) &\mathcal{M}_{\rho ,m}[e^{-ax}g(x)](u\omega,v,\omega)du = \int_0^\infty f(\omega u) \int_0^\infty \frac{e^{-u\omega x-v/x}}{\left(x^m +\omega^m\right)^\rho} e^{-a\omega x} g(\omega x)dx du\\
	&=\int_0^\infty \frac{e^{-a\omega x-v/x}}{\left(x^m +\omega^m\right)^\rho}g(\omega x) \int_0^\infty e^{-u\omega x} f(\omega u)du dx\\
	&=\int_0^\infty \frac{e^{-a\omega x-v/x}}{\left(x^m +\omega^m\right)^\rho}g(\omega x)
	\mathcal{N}[f(u)](\omega x,\omega)dx\\
	&=
	\mathcal{M}_{\rho,m}\big[g(x)\mathcal{N}[f(u)](x,\omega)\big]
	(a\omega,v,\omega).
	\end{align*}
\end{proof}
\begin{corollary}
	Let $f$ and $g$ satisfy the condition \eqref{existencecodition}. Then,
	\begin{align}\label{M-L}
	\int_0^\infty f(u) \mathcal{M}_{\rho ,m}[e^{-ax}g(x)](u\omega,v,\omega)du = 
	\mathcal{M}_{\rho,m}\big[g(x)\mathcal{L}[f(u)](x)\big]
	(a\omega,v,\omega),\, a>0,
	\end{align}
	where $\mathcal{L}$ is the Laplace transform defined by \eqref{Laplace_T}.
\end{corollary}
\begin{corollary}
	Let $f$ satisfy the condition \eqref{existencecodition}. Then,
	\begin{align}\label{M_M}
	\mathfrak{M}\big[\mathcal{M}_{\rho ,m}[e^{-ax}f(x)](u\omega,v,\omega),u\to z\big] = \Gamma(z)
	\mathcal{M}_{\rho,m}\big[x^{-z}f(x)\big] (a\omega,v,\omega),\, a>0,
	\end{align}
	where $\mathfrak{M}$ is the Mellin transform defined by \eqref{Mellin-IT}.
\end{corollary}
\begin{proof}
	It is clear that
	\begin{align*}
	\mathfrak{M}\big[\mathcal{M}_{\rho ,m}[e^{-ax}f(x)](u\omega,v,\omega),u\to z\big] &= \int_0^\infty u^{z-1} \mathcal{M}_{\rho ,m}[e^{-ax}f(x)](u\omega,v,\omega)du \\
	&= \int_0^\infty u^{z-1} \int_0^\infty \frac{e^{-u\omega x-v/x}}{\left(x^m +\omega^m\right)^\rho} e^{-a\omega x} f(\omega x)dx du\\
	&=\int_0^\infty \frac{e^{-a\omega x-v/x}}{\left(x^m +\omega^m\right)^\rho}f(\omega x) \int_0^\infty u^{z-1} e^{-u\omega x}du dx\\
	&=\Gamma(z) \int_0^\infty \frac{e^{-a\omega x-v/x}}{\left(x^m +\omega^m\right)^\rho} \frac{f(\omega x)}{(\omega x)^z}dx\\
	&= \Gamma(z)
	\mathcal{M}_{\rho,m}\left[\frac{f(x)}{x^z}\right]
	(a\omega,v,\omega).
	\end{align*}
\end{proof}
\begin{corollary}
	Let $g$ satisfy the condition \eqref{existencecodition}. Then, for the positive real numbers $a,\nu$ and $\mu$ we have
	\begin{align}\label{B_M}
	\mathfrak{B}_{\nu,\mu}\big[\mathcal{M}_{\rho ,m}[e^{-ax}g(x)]&(u\omega,v,\omega)\big](a\omega,v,\omega)\nonumber\\
	& = \omega^{\nu\mu-1}
	\mathcal{M}_{\rho,m}\left\{\frac{g(x)}{x^{\nu\mu}} H^{1,1}_{1,1} \left[\frac{\omega}{x} \left\vert {\begin{array}{l} (1-\nu\mu,1)\\\left(0,\frac{1}{\nu}\right)\end{array}} \right.\right]\right\} (\omega,v,\omega),
	\end{align}
	where $\mathfrak{B}$ is the Borel-D\v zrbashjan transform defined by \eqref{BorelIT}.
\end{corollary}
\begin{proof}
	From \eqref{newIT} and \eqref{BorelIT}, one gets
	\begin{align}\label{B_M2}
	\mathfrak{B}_{\nu,\mu}\big[\mathcal{M}_{\rho ,m}[e^{-ax}g(x)]&(u\omega,v,\omega)\big](a\omega,v,\omega)\nonumber\\
	&=\nu \omega^{\nu\mu-1} \int_0^\infty e^{-\omega^\nu u^{\nu}} u^{\nu\mu-1} \mathcal{M}_{\rho ,m}[e^{-ax}g(x)](u\omega,v,\omega)du.
	\end{align}
	Operating \eqref{M-N} on the function $f(x)=\nu x^{\nu\mu-1} e^{-x^\nu}$, then 
	\begin{align}\label{B_M3}
	\nu \omega^{\nu\mu-1} \int_0^\infty& e^{-\omega^\nu u^{\nu}} u^{\nu\mu-1} \mathcal{M}_{\rho ,m}[e^{-ax}g(x)](u\omega,v,\omega)du\nonumber\\
	&=\mathcal{M}_{\rho ,m}\left[g(x)\mathcal{N}\left[\nu e^{-u^\nu}u^{\nu\mu-1}\right](x,\omega)\right](a\omega,v,\omega).
	\end{align}
	Here,
	\begin{align*}
	\mathcal{N}\left[\nu e^{-u^\nu}u^{\nu\mu-1}\right](x,\omega)&= \nu\omega^{\nu\mu-1} \int_0^\infty u^{\nu\mu-1} e^{-u^\nu \omega^\nu} e^{-ux} du\\
	&=\nu\omega^{\nu\mu-1} \mathcal{L}\left[u^{\nu\mu-1} e^{-u^\nu \omega^\nu}\right](x)\\
	&=\frac{\nu\omega^{\nu\mu-1}}{x^{\nu\mu}} H^{1,1}_{1,1}\left[\frac{\omega^\nu}{x^\nu}\left\vert {\begin{array}{l}
		(1-\nu\mu,\nu)\\(0,1)
		\end{array}}\right.\right]\\
	&=\frac{\omega^{\nu\mu-1}}{x^{\nu\mu}} H^{1,1}_{1,1}\left[\frac{\omega}{x}\left\vert {\begin{array}{l}
		(1-\nu\mu,1)\\ \left(0,\frac{1}{\nu}\right)
		\end{array}}\right.\right],
	\end{align*}
	where \cite[p. 39, Eq(2.3.12) and p.47, Eq(2.5.25)]{IT24} are used. Thus, \eqref{B_M2} reduces the required result \eqref{B_M}.
\end{proof}
\begin{corollary}
	Let $f$ satisfy the condition \eqref{existencecodition}. Then, the following identity holds true.
	\begin{align}\label{H_M}
	\bf{H}\Big[&\mathcal{M}_{\rho,m}\left[e^{-ax}f(x)\right](u\omega,v,\omega)\Big](v,\omega)\nonumber\\
	&=\mathcal{M}_{\rho,m}\left[\frac{f(x)}{x} H^{m,n+1}_{r+1,s} \left[\frac{\omega}{x}\left\vert {\begin{array}{c}
		(0,1),(a_i ,\alpha_i )_{1,r}\\ (b_j ,\beta_j )_{1,s}
		\end{array}}\right.\right]\right](a\omega,v,\omega),
	\end{align}
	where $\bf{H}$ is the integral transform defined by \eqref{HIT}.
\end{corollary}
\begin{proof}
	From \eqref{newIT} and \eqref{HIT} we have
	\begin{align}\label{H_M2}
	\bf{H}\Big[&\mathcal{M}_{\rho,m}\left[e^{-ax}f(x)\right](u\omega,v,\omega)\Big](v,\omega)\nonumber\\
	&=\int_0^\infty H^{m,n}_{r,s} \left[ \omega u\left\vert {\begin{array}{c}
		(a_i ,\alpha_i )_{1,r}\\ (b_j ,\beta_j )_{1,s}
		\end{array}}\right.\right] \int_0^\infty \frac{e^{-u\omega x -v/x}}{\left(x^m +\omega^m\right)^\rho} e^{-a\omega x} f(\omega x)dx du\nonumber\\
	&=\int_0^\infty \frac{e^{-a\omega x -v/x}}{\left(x^m +\omega^m\right)^\rho} f(\omega x) \int_0^\infty e^{-u\omega x} H^{m,n}_{r,s} \left[ \omega u\left\vert  {\begin{array}{c}
		(a_i ,\alpha_i )_{1,r}\\ (b_j ,\beta_j )_{1,s}
		\end{array}}\right.\right]du dx\nonumber\\
	&=\mathcal{M}_{\rho,m}\left[f(x)\mathcal{N}\left[ H^{m,n}_{r,s} \left[ u\left\vert {\begin{array}{c}
		(a_i ,\alpha_i )_{1,r}\\ (b_j ,\beta_j )_{1,s}
		\end{array}}\right.\right] \right](x,\omega)\right](a\omega,v,\omega)\nonumber\\
	&=\mathcal{M}_{\rho,m}\left[\frac{f(x)}{\omega}\mathcal{L}\left[ H^{m,n}_{r,s} \left[\omega u\left\vert {\begin{array}{c}
		(a_i ,\alpha_i )_{1,r}\\ (b_j ,\beta_j )_{1,s}
		\end{array}}\right.\right] \right](x)\right](a\omega,v,\omega).
	\end{align}
	According to \cite[p.45, Eq. (2.5.16)]{IT24}, we have
	\begin{align*}
	\mathcal{L}\left[ H^{m,n}_{r,s} \left[ u\left\vert  {\begin{array}{c}
		(a_i ,\alpha_i )_{1,r}\\ (b_j ,\beta_j )_{1,s}
		\end{array}}\right.\right] \right](x)
	=\frac{1}{x} 
	H^{m,n+1}_{r+1,s} \left[ \frac{1}{x}\left\vert  {\begin{array}{c}
		(0,1), (a_i ,\alpha_i )_{1,r}\\ (b_j ,\beta_j )_{1,s}
		\end{array}}\right.\right].
	\end{align*}
	Thus, \eqref{H_M2} can be rewritten as the targeted form \eqref{H_M}.
\end{proof}
Next, we are going to obtain a convolution formula for the $\mathcal{M}$-tranform. Let the functions $f(x)$ and $g(x)$ satisfy condition \eqref{existencecodition}. If $F(u,v,\omega)$ and $G(u,v,\omega)$ are the $\mathcal{M}_{\rho,m}$-images of the functions $f(x)$ and $g(x)$, respectively, then
\begin{align*}
F(u,v,\omega)G(u,v,\omega)=& \int_0^\infty \int_0^\infty \frac{f(\omega t)}{\left(t^m +\omega^m\right)^\rho} \frac{g(\omega x)}{\left(x^m +\omega^m\right)^\rho} e^{-u(t+x)} e^{-\frac{v}{t}-\frac{v}{x}} dtdx\\
=& \int_0^\infty \int_x^\infty \frac{f(\omega(y-x))}{\left((y-x)^m +\omega^m\right)^\rho} \frac{g(\omega x)}{\left(x^m +\omega^m\right)^\rho} e^{-uy} e^{-\frac{v}{y-x}-\frac{v}{x}} dydx,
\end{align*}
here, the substitution $y=x+t$ is used. Moreover, by changing the order of the integration, the last double integral equals
\begin{align*}
& \int_0^\infty e^{-uy}\int_0^y \frac{f(\omega(y-x))e^{-\frac{v}{y-x}}} {\left((y-x)^m +\omega^m\right)^\rho} \frac{g(\omega x)e^{-\frac{v}{x}}} {\left(x^m +\omega^m\right)^\rho} dxdy\\
=& \int_0^\infty \frac{e^{-uy-v/y}}{\left(y^m +\omega^m\right)^\rho}\left[\left(y^m +\omega^m\right)^\rho e^{v/y} \int_0^y \frac{f(\omega(y-x))e^{-\frac{v}{y-x}}} {\left((y-x)^m +\omega^m\right)^\rho} \frac{g(\omega x)e^{-\frac{v}{x}}} {\left(x^m +\omega^m\right)^\rho} dx \right]dy.
\end{align*}
Thus, we can define the following $\mathcal{M}_{\rho,m}$-convolution.
\begin{align}\label{M_convolutionFormula}
\left(f*_{_{\mathcal{M}}} g\right)(x)= &\left[\left(\frac{x}{\omega}\right)^m +\omega^m\right]^\rho e^{v\omega/x}\nonumber\\
&\times \int_0^{x/\omega} \frac{f(x-\omega t)e^{-\frac{v\omega}{x-\omega t}}} {\left[\left(\frac{1}{\omega}(x-\omega t)\right)^m +\omega^m\right]^\rho} \frac{g(\omega t)e^{-\frac{v}{t}}} {\left(t^m +\omega^m\right)^\rho} dt.
\end{align}
It is clear that $\left(f*_{_{\mathcal{M}}} g\right)(x)=\left(g*_{_{\mathcal{M}}} f\right)(x)$. Furthermore, the following result holds true.
\begin{theorem}\label{M_convolutionTheorem}
	Let the functions $f(x)$ and $g(x)$ satisfy condition \eqref{existencecodition}. If $F(u,v,\omega)$ and $G(u,v,\omega)$ are the $\mathcal{M}_{\rho,m}$-images of the functions $f(x)$ and $g(x)$, respectively, then
	\begin{align}\label{M_convolutionIdentity}
	F(u,v,\omega)G(u,v,\omega)= \mathcal{M}_{\rho,m}[\left(f*_{_{\mathcal{M}}} g\right)(x)](u,v,\omega).
	\end{align}	
\end{theorem}
%
Next, an inversion formula for the $\mathcal{M}$-transform \eqref{newIT} will be obtained by using the duality relation \eqref{M_L duality}. The inverse Laplace transform is defined (See,\cite{IT8,IT9}) as
\begin{align}
\mathcal{L}^{-1}[F(s)](x)=\frac{1}{2 \pi i} \int_{\alpha-i\infty}^{\alpha+i\infty} e^{sx}F(s)ds,\; \alpha>0.
\end{align}
\begin{theorem}[Inversion Formula]\label{Inversion}
	The inversion of the $\mathcal{M}$-transform \eqref{newIT} is given by
	\begin{align}\label{inversionFormula}
	f(x)=\left(\frac{x^m}{\omega^m}+\omega^m\right)^\rho e^{\omega v/t} \mathcal{L}^{-1}\big[\mathcal{M}_{\rho,m}[f(t)](u,v,\omega)\big]\left(\frac{x}{\omega},\omega,v\right),
	\end{align}
	provided that the involved integrals absolutely converge.
\end{theorem}
\begin{proof}
	Let
	\begin{align}\label{F}
	F(x;v,\omega)=\frac{f(\omega x) e^{-v/x}}{\left(x^m +\omega^m\right)^\rho},\, x>0,
	\end{align}
	where $v\in\mathbb C$, with $\operatorname{Re} v>0$, and $\omega>0$.\\ It is clear that under the condition \eqref{existencecodition}, $F$ is well-defined function. 
	Substituting \eqref{F}  in \eqref{newIT}, gives
	\begin{align*}
	\mathcal{M}_{\rho,m}[f(x)](u,v,\omega)=\mathcal{L}[F(x,v,\omega)](u).
	\end{align*}
	Here, the variables $v,\omega$ of the $\mathcal{M}$-transform are considered as parameters.\\
	Operating the Laplace-inversion formula, yields
	\begin{align*}
	\mathcal{L}^{-1}\left[\mathcal{M}_{\rho,m}[f(x)](u,v,\omega)\right]= F(x,v,\omega).
	\end{align*}
	By replacing $x$ with $x/\omega$ and using \eqref{F}, the required result \eqref{inversionFormula} follows directly.
\end{proof}
\section{Applications}
Here, we discuss two illustrative examples to show the applicability of the new $\mathcal{M}$-transform. 
\begin{example}[First order initial boundary value problem]
	\begin{align}
	w_t +w_x &=p(t,\omega) e^{-v\omega/t}r(t,x),\; t>0,\, x>0,\label{4.1}\\
	w(0,x,\omega)&=\omega\phi(\omega),\; x\geq 0,\label{4.2}\\
	w(t,0,\omega)&=0,\; t\geq 0,\label{4.3}
	\end{align}
\end{example}
where
$$0<p(t,\omega):=\left(\frac{t^m}{\omega^m}+\omega^m\right)^{-\rho},\;\omega>0,\,t\geq 0,\, m\in\mathbb N_0 =\{0,1,2,\cdots\},$$
%
and \(r(t,x)\) and \(\phi(\omega)\) are given functions. Here, \(v\) and \(\omega\) are taken as parameters.\\
%
Equation \eqref{4.1} can be rewritten as
\begin{align}\label{4.4}
\frac{e^{v\omega/t}w_t}{p(t,\omega)} + \frac{e^{v\omega/t}w_x}{p(t,\omega)} =r(t,x),\; t>0,\, x>0.
\end{align}
Applying the $\mathcal{M}$-transform \eqref{newIT} to \eqref{4.4}, and using the duality relation \eqref{elimintateToNaturalIT} of the $\mathcal{M}$-transform with the natural transform, give
\begin{align}\label{4.5}
\frac{u}{\omega} \hat{w}(u,x,v,\omega)-\frac{1}{\omega}w(0,x,\omega)+\hat{w}_x = \mathcal{M}_{\rho,m}[r(t,x)](u,x,v,\omega),
\end{align}
where \(\hat{w}\) is the natural transform image of \(w\). 
Using the initial condition \eqref{4.2}, yields
\begin{align}\label{4.6}
\hat{w}_x +\frac{u}{\omega} \hat{w}=F(x),
\end{align}
where
\begin{align}\label{4.7}
F(x)=\mathcal{M}_{\rho,m}[r(t,x)](u,x,v,\omega)+\phi(\omega).
\end{align}
Transforming the boundary condition \eqref{4.3}, gives
\begin{align}\label{4.8}
\hat{w}(u,0,v,\omega)=0,\, u\geq0.
\end{align}
The solution to problem \eqref{4.6}, \eqref{4.8} is
\begin{align*}
\hat{w}(u,x,v,\omega)=e^{-\frac{u}{\omega}x}\int_0^x F(y) e^{\frac{u}{\omega}y}dy.
\end{align*}
Now, applying the inverse natural transform and using \eqref{4.7}, give
\begin{align}\label{4.10}
w(t,x,v,\omega)=&\mathcal{N}^{-1}\left[e^{-\frac{u}{\omega}x}\int_0^x  e^{\frac{u}{\omega}y} \mathcal{M}_{\rho,m}[r(t,x)](u,y,v,\omega)dy\right]\nonumber\\
& 
+\mathcal{N}^{-1}\left[e^{-\frac{u}{\omega}x} \phi(\omega) \int_0^x e^{\frac{u}{\omega}y}dy\right]\nonumber\\
=&\mathcal{N}^{-1}\left[e^{-\frac{u}{\omega}x}\int_0^x  e^{\frac{u}{\omega}y} \mathcal{M}_{\rho,m}[r(t,x)](u,y,v,\omega)dy\right] \nonumber\\
&
+ \omega \phi(\omega) \mathcal{N}^{-1}\left[\frac{1}{u}\left(1-e^{-\frac{u}{\omega}x}\right) \right] \nonumber\\
=&\mathcal{N}^{-1}\left[e^{-\frac{u}{\omega}x}\int_0^x  e^{\frac{u}{\omega}y} \mathcal{M}_{\rho,m}[r(t,x)](u,y,v,\omega)dy\right] \nonumber\\
&
+ \omega \phi(\omega) \left[\uptheta(t)-\uptheta(t-x) \right],
\end{align}
where $\uptheta(t)$ is the Heaviside function defined as
$$\uptheta(t)=\begin{cases}
1;\quad t>0,\\0;\quad t<0.
\end{cases}$$
\begin{remark}
	In case, \(v=0\) the formula \eqref{4.10} recovers the solution to problem \eqref{4.1}-\eqref{4.3} obtained via the Srivastava-Luo-Raina transform in \cite{IT26}. 
\end{remark}
\begin{example}[Heat flow with sources and homogeneous boundary conditions]\end{example}
\noindent Here we consider the following initial boundary value problem
\begin{equation}\label{2.1-2}
\begin{split}
&\frac{\partial \varphi}{\partial t}= \frac{\partial^{2} \varphi}{\partial x^{2}}+(t^m +1)^{-\rho}e^{-v/t}r(x,t), \quad x \in (0,\pi),\ \ t>0,
\\
&\varphi(x,0)=f(x),\quad x \in [0,\pi],
\\
&\varphi(0,t)=0.\qquad \varphi(\pi,t)=0, \quad t \geq 0, 
\end{split}
\end{equation}
where \(f(x)\) and \(r(x,t)\) are known functions.

Rewriting the differential equation in \eqref{2.1-2} as
\begin{align}\label{2.1-3}
(t^m +1)^{\rho}e^{v/t}\frac{\partial \varphi}{\partial t}-(t^m +1)^{\rho}e^{v/t} \frac{\partial^{2} \varphi}{\partial x^{2}}=r(x,t), \quad x \in (0,\pi),\ \ t>0.
\end{align}
Applying the $\mathcal{M}_{\rho,m}$-transform \eqref{newIT} to \eqref{2.1-3} and using the duality relation \eqref{M_L duality2} with the Laplace transform (when $\omega=1$), give
\begin{align}\label{2.1-4}
\hat{\varphi}_{xx}-u \hat{\varphi}(x,u)=-\mathcal{M}_{\rho,m}[r(x,t)](u,v,1)-f(x),\; x\in[0,\pi],
\end{align}
where \(\hat{\varphi}\) is the Laplace transform of \(\varphi\). Transforming the boundary data in \eqref{2.1-2}, gives
\begin{align}\label{2.1-5}
\hat{\varphi}(0,u)=\hat{\varphi}(\pi,u)=0.
\end{align}
Using the variation of parameters method, one can obtain the solution to the boundary value problem defined by \eqref{2.1-4} and \eqref{2.1-5} as
\begin{align}\label{hat_phi}
\hat{\varphi}(x,u)=&-\frac{1}{\sqrt{u}}\int_0^x \frac{\sinh\sqrt{u}x \sinh\sqrt{u}(\pi-y)-\sinh\sqrt{u}\pi \sinh\sqrt{u}(x-y)}{\sinh\sqrt{u}\pi}G(y)dy\nonumber\\
&-\frac{1}{\sqrt{u}}\int_x^\pi \frac{\sinh\sqrt{u}x \sinh\sqrt{u}(\pi-y)}{\sinh\sqrt{u}\pi}G(y)dy,
\end{align}
where
\begin{align}\label{G}
G(x)&=-\mathcal{M}_{\rho,m}[r(x,t)]-f(x)\nonumber\\
&=-\mathcal{L}\left[\frac{e^{-v/t}r(x,t)}{\left(t^m +1\right)^{\rho}}\right](u)-f(x),\; x\in[0,\pi].
\end{align}
Here, the daulity relation \eqref{M_L duality} is used.

Since,
\[\Delta(u):=\frac{\sinh\sqrt{u}x \sinh\sqrt{u}(\pi-y)}{\sqrt{u}\sinh\sqrt{u}\pi}=O(1)\; as\; u\to 0,\]
then, \(\Delta(u)\) and hence \(\hat{\varphi}(x,u)\) have simple poles at \(u_k = -k^2\) for all \(k=1,2,\cdots\).
Therefore,
\begin{align}\label{L_Theta}
\mathcal{L}^{-1}\left[\frac{\sinh\sqrt{u}x \sinh\sqrt{u}(\pi-y)}{\sqrt{u}\sinh\sqrt{u}\pi}\right]=\Theta(x,t;y),
\end{align}
with
\begin{align}\label{Theta}
\Theta(x,t;y)=\frac{2}{\pi} \sum_{k=1}^{\infty} (-1)^{k+1} e^{-k^2 t}\sin kx \sin ky.
\end{align}
In view of \eqref{G}, \eqref{L_Theta} and \eqref{Theta}, \eqref{hat_phi} can be rewritten as
\begin{align}\label{hat_phi2}
\hat{\varphi}(x,u)=-\int_0^\pi \mathcal{L}[\Theta(x,t;y)] \mathcal{L}\left[\frac{e^{-v/t}r(y,t)}{\left(t^m +1\right)^{\rho}}\right]dy
-\int_0^\pi \mathcal{L}[\Theta(x,t;y)] f(y)dy.
\end{align}
By the Laplace convolution theorem, the inverse Laplace transform of \eqref{hat_phi2} is
\begin{align}\label{phi}
\varphi(x,t)=&-\int_0^\pi \int_0^t \Theta(x,t-\zeta;y) \frac{e^{-v/\zeta}r(y,\zeta)}{\left(\zeta^m +1\right)^{\rho}}d\zeta dy
-\int_0^\pi \Theta(x,t;y) f(y)dy.
\end{align}
Substituting from \eqref{Theta} into \eqref{phi}, gives 
\begin{align}\label{phi-final}
\varphi(x,t)=&\frac{2}{\pi} \sum_{k=1}^{\infty} (-1)^{k}\left(\int_0^\pi \int_0^t    \frac{e^{k^2 \zeta-v/\zeta}r(y,\zeta)}{\left(\zeta^m +1\right)^{\rho}}\sin kyd\zeta dy\right)e^{-k^2 t}\sin kx
\nonumber\\
&+\frac{2}{\pi} \sum_{k=1}^{\infty} (-1)^{k} \left(\int_0^{\pi}\sin ky f(y)dy\right)e^{-k^2 t}\sin kx,
\end{align}
that is the solution to Problem \eqref{2.1-2}. 

In particular, if
%
\(r(x,t)=e^{-t}\sin 3x, \; x\in[0,\pi],\, t\geq0,\)
we find
\begin{align*}
\varphi(x,t)=&\frac{2}{\pi} \sum_{k=1}^{\infty} (-1)^{k}\left(\int_0^\pi \int_0^t    \frac{e^{(k^2 -1) \zeta-v/\zeta}}{\left(\zeta^m +1\right)^{\rho}}\sin 3y \sin ky d\zeta dy\right)e^{-k^2 t}\sin kx
\nonumber\\
&+\frac{2}{\pi} \sum_{k=1}^{\infty} (-1)^{k} \left(\int_0^{\pi}\sin ky f(y)dy\right)e^{-k^2 t}\sin kx.
\end{align*}
Since,
\[\int_0^\pi \sin 3y \sin ky dy=
\begin{cases}
\frac{\pi}{2},& k=3\\ 0,& k\neq 3,
\end{cases}\]
then
\begin{align}\label{phi-special}
\varphi(x,t)=&-\left( \int_0^t    \frac{e^{8\zeta-v/\zeta}}{\left(\zeta^m +1\right)^{\rho}}d\zeta\right)e^{-9 t}\sin kx
\nonumber\\
&+\frac{2}{\pi} \sum_{k=1}^{\infty} (-1)^{k} \left(\int_0^{\pi}\sin ky f(y)dy\right)e^{-k^2 t}\sin kx.
\end{align}
%
\begin{remark}
	$\bullet$
	Formula \eqref{phi-special} recovers the solution to Problem \eqref{2.1-2}, when $v=0$ and $r(x,t)=e^{-t}\sin 3x$, obtained via the Srivastava-Luo-Raina integral transform method in \cite{IT1}.
	
	$\bullet$
	Formula \eqref{phi-special} recovers the solution to problem \eqref{2.1-2}, when $\rho=0,\,v=0$ and $r(x,t)=e^{-t}\sin 3x$, obtained via the Eigenfunction expansion method in \cite{IT27}. 
\end{remark}
\section{Conclusion}
In this article a new integral transform is introduced, which depends on a number of parameters  so that it covers many known integral transforms as its special cases. In Section 2, duality relations of the new transform with well-known transforms, such as the Laplace transform, the natural transform and the Srivastava-Luo-Riana transform have been shown, some examples have been discussed. Also, an extended H-function has been defined. 
Integral indentities involving the new transform, well-known integral transforms and special functions are established in Section 3. These identities include Parseval-type indentity, Mellin integrals, $\bf{H}$-functions and the $\bf{H}$-transform. Convlution theorem and inversion theorem for the new transform are also estabilshed. 

The new transform is a precious tool for solving certain initial and boundary value problems with certain variable coefficients. We have shown, in Section 4, that an appropriate choice of the parameters of the new transform helps to  eliminate the variable coefficients in the problem. And hence by using the duality relations of the new transform with well-known transforms, such as the Laplace transform and the natural transform, we can proceed to find the solutions.
%
\section*{Acknowledgment}
The author thanks the editors and the referees for their valuable comments and suggestions which improved greatly the quality of this paper.
\section*{Declarations}

%

\subsection*{Availability of data and material}
Not applicable.

\subsection*{Competing interests} 
The author declares that he has no competing interests.

\subsection*{Funding}Not applicable.

\subsection*{Authors' contributions}
Unique contributor. The author read and approved the final manuscript.

\end{document}